\documentclass{elsart3-1}

\usepackage{amssymb}
\usepackage{amsfonts, amsmath, psfrag, graphicx, subfigure, url}
\usepackage[english,francais]{babel}

\DeclareMathOperator{\grad}{grad}

\newcommand{\fbnod}[1]{\phi^h_{#1}}
\newcommand{\fbedg}[1]{\lambda^h_{#1}}
\newcommand{\cbnod}[1]{\phi^H_{#1}}
\newcommand{\cbedg}[1]{\lambda^H_{#1}}
\newcommand{\card}[1]{\lvert #1 \lvert}

\newcommand{\tL}{\tilde{L}}
\newcommand{\tC}{\tilde{C}}
\newcommand{\bC}{\hat{C}}
\newcommand{\tI}{\tilde{I}}
\newcommand{\tF}{\tilde{F}}

\newcommand{\Srond}{\mathcal{S}}
\newcommand{\Trond}{\mathcal{T}}
\newcommand{\Erond}{\mathcal{E}}

\renewcommand{\tilde}{\widetilde}
\renewcommand{\hat}{\widehat}

\newtheorem{theorem}{Theorem}[section]
\newtheorem{lemma}[theorem]{Lemma}
\newtheorem{e-proposition}[theorem]{Proposition}

\newtheorem{e-definition}[theorem]{Definition\rm}

\renewcommand{\theequation}{\arabic{equation}}
\def\og{\leavevmode\raise.3ex\hbox{$\scriptscriptstyle\langle\!\langle$~}}
\def\fg{\leavevmode\raise.3ex\hbox{~$\!\scriptscriptstyle\,\rangle\!\rangle$}}

\begin{document}

\centerline{Numerical Analysis
}
\begin{frontmatter}



\selectlanguage{english}
\title{Gradient-prolongation commutativity and graph theory}

\author[label1]{Fran{
    c}ois Musy},
\ead{francois.musy@ec-lyon.fr}
\author[label2]{Laurent Nicolas},
\ead{laurent.nicolas@ec-lyon.fr}
\author[label1,label2]{Ronan Perrussel}
\ead{ronan.perrussel@ec-lyon.fr}
\address[label1]{Institut Camille Jordan, Ecole Centrale de Lyon,
F-69134 Ecully Cedex} 
\address[label2]{Centre de G\'enie \'Electrique de Lyon, Ecole
Centrale de Lyon, F-69134 Ecully Cedex} 

 
\medskip
\begin{center}
  {\small Received *****; accepted after revision +++++\\
    Presented by £££££}
\end{center}
 
\begin{abstract}
  This note gives conditions that must be imposed to algebraic multilevel
  discretizations involving at the same time nodal and edge elements so that a
  gradient-prolongation commutativity condition will be satisfied; this
  condition is very important, since it characterizes the gradients of coarse
  nodal functions in the coarse edge function space.  They will be expressed
  using graph theory and they provide techniques to compute approximation
  bases at each level.  {\it To cite this article: A.  Name1, A.  Name2, C. R.
    Acad.  Sci. Paris, Ser. I 340 (2005).}
  
  \vskip 0.5\baselineskip
  
  \selectlanguage{francais}
\noindent{\bf R\'esum\'e}
\vskip 0.5\baselineskip
\noindent
{\bf Commutativit\'e entre gradient et prolongement et th\'eorie des graphes}
Cette note donne des conditions qui doivent \^etre impos\'ees aux
discr\'etisations multiniveau alg\'ebriques en \'el\'ements finis nodaux et
d'ar\^ete de fa\c con \`a assurer la commutativit\'e entre gradient et
prolongement; cette relation importante caract\'erise les gradients des
fonctions nodales grossi\`eres dans l'espace des fonctions d'ar\^ete
grossi\`eres.  Ces conditions seront exprim\'ees en terme de graphes et elles
permettent d'introduire des m\'ethodes de calcul des bases d'approximation aux
diff\'erents niveaux.  {\it Pour citer cet article~: A.  Name1, A. Name2, C.
  R.  Acad. Sci. Paris, Ser.  I 340 (2005).}

\end{abstract}

\end{frontmatter}

\selectlanguage{francais}
\section*{Version fran\c{c}aise abr\'eg\'ee}

L'approximation num\'erique du champ \'electrique ou magn\'etique utilise
fr\'e\-quem\-ment les \'el\'ements finis d'ar\^ete dont la relation avec les
\'el\'ements finis nodaux traduit des propri\'et\'es importantes au niveau
discret~\cite{HIP02}.  Dans ce qui suit, nous consid\`ererons les \'el\'ements
de plus bas degr\'e~: $P_1$ en nodal et ordre $1$ incomplet pour les ar\^etes.
D\`es qu'on traite des probl\`emes de grande taille, une strat\'egie
multiniveau est un choix int\'eressant.  Pour les syst\`emes provenant de
discr\'etisations par \'el\'ements finis d'ar\^ete, Hiptmair a introduit des
m\'ethodes multiniveau pour une hi\'erarchie de maillages
embo\^it\'es~\cite{HIP99}.

Cependant, dans des applications r\'ealistes, on ne dispose g\'en\'eralement
pas de maillages structur\'es.  La strat\'egie multiniveau alg\'ebrique va
donc s'imposer~: il s'agit de d\'efinir des fonctions grossi\`eres nodales et
d'ar\^ete gr\^ace aux contributions de paquets de fonctions fines nodales et
d'ar\^ete; les combinaisons lin\'eaires \eqref{eq:alpha_unk} et
\eqref{eq:beta_unk} d\'efinissent respectivement ces fonctions grossi\`eres
nodales et d'ar\^ete.

Par construction les gradients des fonctions nodales fines appartiennent \`a
l'espace des fonctions d'ar\^ete fines ce que traduit la relation
\eqref{eq:grad_fin}.  Dans cette relation, $G^h$ est la matrice d'incidence
arcs-sommets du graphe orient\'e naturellement associ\'e au maillage de
travail. 
Les orientations des {arcs} sont arbitraires.

Pour adapter aux m\'ethodes alg\'ebriques les lisseurs des m\'ethodes
g\'eom\'etriques, Reitzinger et Sch\"oberl~\cite{REI02} ont introduit une
repr\'esentation explicite des gradients des fonctions grossi\`eres nodales
dans la base des fonctions grossi\`eres d'ar\^ete, donn\'ee par la
relation~\eqref{eq:grad_gros} o\`u $G^H$ est une matrice d'incidence
arcs-sommets.

En regroupant les relations \eqref{eq:prolongement} \`a \eqref{eq:grad_gros},
nous obtenons la relation matricielle \eqref{eq:matrices_compatibles}.  La
matrice $\alpha$ est construite par exemple par les m\'ethodes d\'efinies dans
\cite{MAN99} qui permettent d'obtenir les fonctions grossi\`eres nodales comme
partition de l'unit\'e et de contraindre leurs supports \`a \^etre inclus dans
des ensembles g\'eom\'etriques convenablement choisis.

Connaissant $G^h$ et $\alpha$, nous souhaitons choisir $G^H$ comme matrice
d'incidence arcs-sommets d'un graphe orient\'e $\Srond^H$. Nous donnons dans
cette note une condition n\'ecessaire et suffisante sur ce graphe, la
proposition \eqref{prop:principale}, qui assure l'existence d'une solution de
\eqref{eq:matrices_compatibles}.  En effet, nous associerons, par un
proc\'ed\'e d\'ecrit dans la partie en anglais, \`a chaque ar\^ete fine un
sous-graphe du graphe grossier, qui doit \^etre connexe.

La connaissance de ces sous-graphes donne les degr\'es de libert\'e
disponibles pour d\'eterminer des fonctions d'ar\^ete grossi\`eres compatibles
avec les fonctions nodales grossi\`eres; en r\'esolvant un probl\`eme de flot
sur ces sous-graphes, voir \eqref{eq:system_betai}, nous pouvons alors
construire la matrice $\beta$ (Section \ref{sec:construit_beta}).

\selectlanguage{english}
\section{Introduction}
\label{sec:introduction}
Numerical approximation of electric or magnetic field uses often edge finite
elements whose relation with nodal {finite} elements contains important
properties at discrete level~\cite{HIP02}.  In this note we restrict ourselves
to lowest order approximation~: $P_1$ for nodal elements and incomplete order
$1$ for edge elements.  In order to solve large problems, multilevel methods
are an attractive choice.  While, for systems coming from edge element
discretisation, Hiptmair~\cite{HIP99} proposed multilevel methods using nested
meshes, engineering applications do not usually provide structured meshes.
Therefore, algebraic multilevel methods are an interesting option: we have to
build coarse nodal and edge functions by using aggregates of fine nodal and
edge functions. If $(\fbnod{p})_{p=1,\ldots,N^h}$ and
$(\fbedg{i})_{i=1,\ldots,E^h}$ respectively denote fine nodal and edge bases,
the following linear combinations define coarse nodal and edge functions:
\begin{subequations}
  \begin{align}
    \cbnod{n} & = \sum_{p=1}^{N^h} \alpha_{pn} \fbnod{p}, \ \forall n
    \in \{1, \ldots N^H\}, \label{eq:alpha_unk} \\
    \cbedg{e} & = \sum_{i=1}^{E^h} \beta_{ie} \fbedg{i}, \ \forall e \in \{1,
    \ldots, E^H \}. \label{eq:beta_unk}
  \end{align}
  \label{eq:prolongement}
\end{subequations}

By construction, the gradients of fine nodal functions belong to the space of
fine edge functions:
\begin{equation}
  \forall p \in \{1,\ldots,N^h\}, \ \grad(\fbnod{p}) = \sum_{i=1}^{E^h} G^h_{ip} \fbedg{i},
  \label{eq:grad_fin}
\end{equation}
where $G^h$ is the edge-node incidence matrix of the digraph naturally
associated {with} the initial mesh.  The orientation of the edges can be
arbitrarily chosen.

In \cite{REI02}, Reitzinger and Sch\"oberl deduced their smoother from the
matrix $G^H$ involved in the relation:
\begin{equation}
  \forall n \in \{1,\ldots,N^H\}, \ \grad(\cbnod{n}) = \sum_{e=1}^{E^H} G^H_{en} \cbedg{e},
  \label{eq:grad_gros}
\end{equation}
which states that the gradients of the coarse nodal functions must also belong
to the space of coarse edge functions.  The matrix $G^H$ is an edge-node
incidence matrix as in the structured case.  Relation~\eqref{eq:grad_gros}
does not guarantee the efficacy of the algebraic multilevel method but it
leads to relevant strategies.

Gathering Equations \eqref{eq:prolongement}, \eqref{eq:grad_fin} and
\eqref{eq:grad_gros}, we obtain the matrix relation:
\begin{equation}
  G^h \alpha = \beta G^H.
  \label{eq:matrices_compatibles}
\end{equation}
The matrix $\alpha$ is constructed following {for instance} the methods
defined in \cite{MAN99}, which provides a family of coarse nodal functions,
making up a partition of unity, whose supports satisfy appropriate conditions.

Knowing the left-hand side of \eqref{eq:matrices_compatibles}, we want to
choose $G^H$ as an edge-node incidence matrix of a digraph $\Srond^H$, and we
will give conditions on the coarse graph $\Srond^H$, which ensure the
existence of a matrix $\beta$ satisfying \eqref{eq:matrices_compatibles}.
Moreover, the proof of the proposition indicates how to choose the degrees of
freedom which enables us to define the coarse edge functions. It also helps us
to construct \(\beta\).

\section{Notation and statement of the problem}

Let $(L_n)_{n=1,\ldots,N^H}$ be sets of indices in $\{1,\ldots,N^h\}$ such
that:
\begin{equation}
   \bigcup_{n=1}^{N^H} L_n = \{1,\ldots,N^h\}. \label{eq:condition_Ln}
\end{equation}

\begin{figure}[htbp]
  \centering \subfigure[Fine graph with node indices. Sets
  $(L_n)_{n=1,\ldots,3}$ are surrounded.]{
    \label{fig:Ln}
    \includegraphics[width=0.25\textwidth]{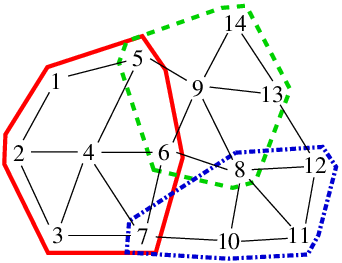}
  }
  \subfigure[Fine graph with edge indices. Set $C_3$ is represented by bold edges.]{
    \label{fig:C3}
    \includegraphics[width=0.25\textwidth]{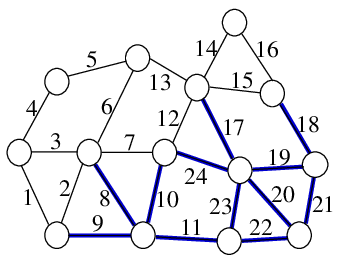}
  }
  \subfigure[Coarse graph.]{
    \label{fig:graphegros}
    \includegraphics[width=0.17\textwidth]{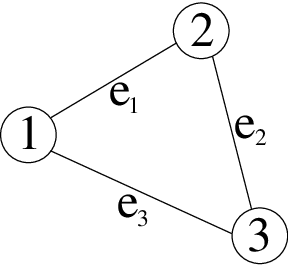}
  }
  \subfigure[Set $I_{e_3}$ is represented by bold edges.]{
    \label{fig:Ie}
    \includegraphics[width=0.25\textwidth]{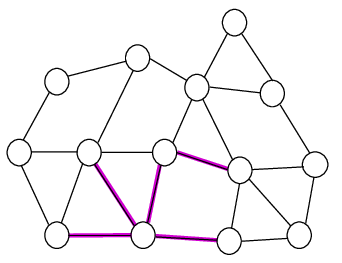}
  }
  \caption{Representation of the fine and coarse graphes, sets $(L_n)_{n=1,\ldots,3}$, $C_3$ and $I_{e_3}$.}
\end{figure}

The matrix $\alpha$ describes the coarse nodal basis; we assume that
{it} is has been previously computed and it has the following
properties:
\begin{itemize}
\item the coarse nodal functions make up a partition of unity, which can be
  algebraically stated as:
  \begin{equation}
    \forall p \in \{1,\ldots,N^h\}, \ \sum_{n=1}^{N^H} \alpha_{pn} = 1,
    \label{eq:partition_unite}
  \end{equation}
\item in order to restrict the support of each coarse basis function
  $\cbnod{n}$, the indices of the non-zero components of $\cbnod{n}$ are
  included in the set $L_n$, \textsl{i.e.}:
  \begin{equation}
    p \in \{1,\ldots,N^h\} \setminus L_n \implies \alpha_{pn} = 0.
    \label{eq:local_support}
  \end{equation}
  The fine nodal function $\fbnod{p}$ {contribute}s to the coarse nodal
  function $\cbnod{n}$ if $p$ belongs to $L_n$.
\end{itemize}
We have a reciprocal set-valued function $\tL$: the set $\tL_p$ is the set of
coarse nodal function indices to which the fine nodal function $\fbnod{p}$
contributes.  For the fine graph in Figure \ref{fig:Ln}, we set $L_1 = \{1, 2,
3, 4, 5, 6, 7\}$, $L_2 = \{5, 6, 8, 9, 13, 14\}$ and $L_3 = \{7, 8, 10, 11,
12\}$. One obtains, for instance, the set $\tL_7 = \{1, 3\}$.

We define two families of sets of fine edge function indices.  We will denote a
directed fine edge $i$ by $\overline{pq}^h$ where $p$ and $q$ are respectively
the starting and ending nodes of the edge $i$.  A similar notation is used for
a directed coarse edge $e=\overline{mn}^H$.

The set $C_n$ is the set of indices of fine edges which have an extremity in
$L_n$:
\begin{equation}
  \label{eq:def_Cn}
  C_n = \bigl\{ i \in \{1,\ldots,E^h\}: i = \overline{pq}^h, \ p \in L_n
  \text{ or } q \in L_n \bigr\}.
\end{equation}

The fine edge function $\fbedg{i}$ {contribute}s to the gradient of the
coarse nodal function $\cbnod{n}$ if $i$ belongs to $C_n$.  Indeed, for the
directed fine edge $i=\overline{pq}^h$, $G^h_{ir}$ is equal to $-1$ if $r = p$
and $+1$ if $r = q$.  Moreover, if $p$ and $q$ are not in $L_n$, the
components $\alpha_{pn}$ and $\alpha_{qn}$ vanish according to
\eqref{eq:local_support}; therefore:
\begin{equation}
  \label{eq:participe_Cn}
  i \in \{1,\ldots,E^h\}\setminus C_n \implies (G^h \alpha_{\bullet n})_i = 0,
\end{equation}
where $\alpha_{\bullet n}$ denotes the $n$-th column of $\alpha$.  The
reciprocal set-valued function $\tC$ is such that $\tC_i$ is the set of coarse
nodal function indices to whose gradient the fine edge function $\fbedg{i}$
{contribute}s.  On Figure \ref{fig:C3}, the fine edges are numbered, set $C_3$
is highlighted and we can note, for instance, the set $\tC_{8} =\{1, 3\}$.

Let \(e=\overline{mn}^H\) be an edge of the coarse graph \(\Srond^H\); we
define:
\begin{equation}\label{eq:def_Ie}
  I_e = C_n \cap C_m.
\end{equation}
By analogy with the structured case and for restricting the support of
$\cbedg{e}$, we enforce:
\begin{equation}\label{eq:participe_Ie}
  i \in \{1,\ldots,E^h\} \setminus I_e \implies \beta_{ie} = 0.
\end{equation}
The fine edge function $\fbedg{i}$ contributes to the coarse edge function
$\cbedg{e}$ if $i$ belongs to $I_e$.  The set-valued function $\tI$ is such
that $\tI_i$ is the set of coarse edge function indices to which the fine edge
function $\fbedg{i}$ {contribute}s. The coarse graph in
Figure~\ref{fig:graphegros} is related to the fine in Figure~\ref{fig:Ln}.
Set $I_{e_3}$ is represented in Figure \ref{fig:Ie}.

The following statement can be easily deduced from \eqref{eq:def_Cn} and the
definition of \(G^h\):

\begin{lemma} \label{lem:tCi}
  If \(i\) denotes the edge \(\overline{pq}^h\), \(\tilde C_i=\tilde L_p\cup
  \tilde L_q\).
\end{lemma}

In order to simplify notations, we introduce the set $ \tilde F= \bigl\{i\in
\{1,\dots, E^h\}: \tilde I_i\neq \emptyset\bigr\}$, since some fine edge
functions might not {contribute} to any coarse edge functions.

For any fine edge \(i\), let \(\Srond^{H, i}\) be the induced subgraph defined
by \(\tilde C_i\): the vertices of \(\Srond^{H, i}\) are the vertices of
\(\Srond^H\), which are indexed by the elements of \(\tilde C_i\) and the
edges of \(\Srond^{H, i}\) are those edges of \(\Srond^H\) whose extremities
are vertices of \(\Srond^{H, i}\).

The following lemma is a direct consequence of definition
\eqref{eq:def_Ie}: 

\begin{lemma} \label{lem:SHi}
  For any edge \(i\in \tilde F\), the edges of \(\Srond^{H, i}\) are
  those edges of \(\Srond^H\) which are indexed by \(\tilde I_i\).
\end{lemma}

We may now state precisely our main result, which gives a necessary and
sufficient condition on the coarse graph \(\Srond^H\) permitting the
resolution of \eqref{eq:matrices_compatibles}:

\begin{e-proposition}
  For all matrices \(\alpha\) satisfying conditions \eqref{eq:partition_unite}
  and \eqref{eq:local_support}, there exists a matrix \(\beta\) satisfying
  \eqref{eq:participe_Ie} and solving \eqref{eq:matrices_compatibles} iff for
  all \(i\), the induced subgraph \(\Srond^{H, i}\) is connected.
  \label{prop:principale}
\end{e-proposition}

\section{The essential steps of the proof}\label{sec:strategy-proof}

\textbf{First step.}  Many relations in \eqref{eq:matrices_compatibles} reduce
to $0=0$: this is the case for $n \notin \tC_i$.

Indeed, according to \eqref{eq:participe_Cn} the $(i, n)$ coefficient of the
right-hand side of \eqref{eq:matrices_compatibles} vanishes.

Conversely, if $e$ does not belong to $\tI_i$, according to
\eqref{eq:participe_Ie} and the definition of $\tI_i$, $\beta_{ie}$ vanishes
and:
\begin{equation}
  \label{eq:reduit_sombeta}
  \sum_{e=1}^{E^H} \beta_{ie} G^H_{en} = \sum_{e \in \tI_i} \beta_{ie} G^H_{en}.
\end{equation}
On the other hand if the directed coarse edge $e$ denoted by $\overline{lm}^H$
belongs to $\tI_i$, Lemma~\ref{lem:SHi} implies that $l$ and $m$ belongs to
$\tC_i$.  However, for $G^H_{en}$ not to vanish for all $e$, $m$ or $l$ must
be equal to $n$, which means that $n$ belong to $\tC_i$, and this contradicts
the assumption $n \notin \tC_i$.

\textbf{Second step.}  We look at all the other equations, \textsl{i.e.} those
for which $n \in \tC_i$.  We note that \eqref{eq:reduit_sombeta} remains and
that the edges indexed by $\tI_i$ are precisely those of the graph
$\Srond^{H,i}$ according to Lemma~\ref{lem:SHi}.

We assume now $i \in \tF$ and we define $G^{H, i}$ as the edge-node incidence
matrix of $\Srond^{H, i}$ and the $(i, n)$ equation of
\eqref{eq:matrices_compatibles} is rewritten:
\begin{equation}
  \label{eq:local_systems}
  \sum_{e \in \tI_i} \beta_{ie} G_{en}^{H,i} = \Theta_{i, n} \text{ where }
  \Theta_{i, n} = \sum_{r \in L_n} G^h_{ir} \alpha_{rn}.
\end{equation}
This could be satisfied for all couples $(i, n)$ such that $n \in
\{1,\ldots,N^H\}$ and $i \in C_n$ or equivalently $i \in \{1,\ldots,E^h\}$ and
$n \in \tC_i$. For a fixed $i$, we may write that $\beta_{i \bullet}$, the
$i$-th row of $\beta$ satisfies the system:
\begin{equation}
  \label{eq:system_betai}
  \sum_{e \in \tI_i} \beta_{ie} G^{H,i}_{en} = \Theta_{i, n}, \ \forall n \in \tC_i.
\end{equation}
{Thus}, we solve line by line for $\beta$ and we see that
\eqref{eq:system_betai} is a 
flow problem whose solution is of the form:
\begin{equation}
  \label{eq:decompose_betai}
  \beta_{i\bullet} = \beta^{\prime}_{i\bullet} + \beta^{\prime \prime}_{i\bullet}.  
\end{equation}
with $(\beta^{\prime \prime}_{i\bullet})^t \in \ker(G^{H,i})^t$ and
$\beta_{i\bullet}^{\prime}$ a particular solution.

More precisely, let $\Trond^i$ be a spanning tree for $\Srond^{H, i}$; call
$\Gamma^i$ the edge-node incidence matrix associated with $\Trond^i$; we know
that $\Gamma^i$ has $\card{\tC_i}-1$ rows and $\card{\tC_i}$ columns, and it
is of rank $\card{\tC_i}-1$.  We choose a vertex $m$ in $\Gamma^i$ and we
solve the system:
\begin{equation}
  \label{eq:system_regulier}
  \sum_{e \in \Erond(\Trond^i)} \beta^{\prime}_{ie} \Gamma^i_{en} = \Theta_{i,
  n}, \ \forall n \in \tC_i \setminus \{m\},
\end{equation}
where $\Erond(\Trond^i)$ denotes the set of indices of the edges of
$\Trond^i$.  The system~\eqref{eq:system_regulier} is a regular system of
$\card{\tC_i} - 1$ equations with $\card{\tC_i} - 1$ unknowns, and we put
{$\beta^{\prime}_{ie}$ equal to $0$ if $e$ is in $\tI_i \setminus
  \Erond(\Trond^i)$}.

It remains to show that the forgotten equation of index $m$ in
\eqref{eq:system_regulier} is automatically satisfied.
Indeed, by denoting $i$ by $\overline{pq}^h$, we sum the right-hand side of
\eqref{eq:system_betai} with respect to $n \in \tC_i$:
\begin{equation}
  \label{eq:somme_tCi}
  \sum_{n \in \tC_i} \Theta_{i, n} = \sum_{n \in \tL_p
  \cup \tL_q} \alpha_{qn} - \alpha_{pn} = 0,
\end{equation}
since in view of \eqref{eq:partition_unite} and \eqref{eq:local_support}, $
\displaystyle \sum_{n \in \tL_p \cup \tL_q} \alpha_{pn} = \sum_{n \in \tL_p}
\alpha_{pn} = 1$ and $\displaystyle \sum_{n \in \tL_p \cup \tL_q} \alpha_{qn}
= \sum_{n \in \tL_q} \alpha_{qn} = 1$.

On the other hand, if we sum the left-hand side of \eqref{eq:system_betai}
with respect to $n \in \tC_i$, we obtain:
\begin{equation}
  \label{eq:somme_tCi2}
  \sum_{n \in \tC_i} \sum_{e \in \tI_i} \beta_{ie} G^{H,i}_{en} = \sum_{e \in
  \tI_i} \beta_{ie} \sum_{n \in \tC_i} G^{H, i}_{en} = 0,
\end{equation}
since each line of $G^{H,i}$ contains only two non-zero coefficients $+1$ and
$-1$.

If $i \notin \tF$, $\card{\tC_i} = \card{\tL_p} = \card{\tL_q} = 1$ and the
relation $\displaystyle \sum_{e=1}^{E^H} \beta_{ie} G^H_{en} = \Theta_{i, n}$
is satisfied from \eqref{eq:reduit_sombeta} and \eqref{eq:somme_tCi}.

Now we assume that $\Srond^{H, i}$ is not connected and we denote by $\bC_i$
the nodes of a connected component.  For the same reasons as in
\eqref{eq:somme_tCi2}, if $\beta$ satisfies \eqref{eq:participe_Ie} one gets
$ \displaystyle \sum_{n \in \bC_{i}} \sum_{e=1}^{E^H} \beta_{ie} G^H_{en} =
0.$

However we can construct a matrix $\alpha$ satisfying
\eqref{eq:partition_unite} and \eqref{eq:local_support} such that $
\displaystyle \sum_{n \in \bC_i} \sum_{r=1}^{N^h} G^h_{ir} \alpha_{rn} \neq
0$.  In fact, in view of \eqref{eq:local_support}, for $i = \overline{pq}^h$
we can write:
\begin{equation}
  \label{eq:somme_tCil}
  \sum_{n \in \bC_i} \sum_{r=1}^{N^h} G^h_{ir} \alpha_{rn} = \sum_{n \in \bC_i} \alpha_{qn} - \alpha_{pn} = \sum_{\bC_i \cap
  \tL_q} \alpha_{qn} -  \sum_{\bC_i \cap \tL_p} \alpha_{pn}.
\end{equation}
Since $\bC_i$ is strictly included in $\tL_p \cup \tL_q$, we will have
$\bC_i \cap \tL_p \neq \tL_p \text{ or } \bC_i \cap \tL_q \neq \tL_q$.
Depending on the situation, we can construct a suitable matrix $\alpha$ such
that:
\begin{equation*}
  \bigg( \sum_{\bC_i \cap \tL_q} \alpha_{qn} = 1 \text{ and }
  \sum_{\bC_i \cap \tL_p} \alpha_{pn} = 0 \bigg) \text{ or } \bigg(
  \sum_{\bC_i \cap \tL_q} \alpha_{qn} = 0 \text{ and } \sum_{\bC_i \cap \tL_p} \alpha_{pn} = 1 \bigg).
\end{equation*}
For these matrices $\alpha$, the condition defined by
\eqref{eq:matrices_compatibles} cannot be ensured.

\section{Construction of the coarse edge functions}
\label{sec:construit_beta}

For a coarse graph satisfying the condition of
Proposition~\ref{prop:principale} and by using the decomposition
\eqref{eq:decompose_betai}, any compatible matrice can be written $\beta =
\beta^{\prime} + \beta^{\prime \prime}$, where the complete matrices are
defined by gathering the lines of index $i$ $\beta_{i\bullet}$,
$\beta^{\prime}_{i\bullet}$ and $\beta^{\prime \prime}_{i\bullet}$.  The
computation of each $\beta^{\prime}_{i \bullet}$ can be done by solving system
\eqref{eq:system_regulier}. As concerns $\beta^{\prime \prime}_{i\bullet}$, a
basis of the kernel of $(G^{H,i})^t$ is given by a set of $k_i$
independent cycles of $\Srond^{H,i}$
.
Then, $\sum_{i \in \tF} k_i$ degrees of freedom should be determined by
minimising an appropriate energy functional; such a problem is introduced
in~\cite{MAP394} and can be related to explanations in~\cite{MAN99}.

We thank Michelle Schatzman for many fertile discussions.

\end{document}